\documentclass[a4pitaper,11pt]{amsart}
\usepackage{amsmath,amsthm,amssymb}
\usepackage[mathscr]{eucal}
 \usepackage{cite}
\usepackage{upgreek}
\usepackage[bookmarks,bookmarksnumbered]{hyperref}
\usepackage{enumerate}

\numberwithin{equation}{section}

\newcommand{\car}{\curvearrowright}

\theoremstyle{plain}
\newtheorem{main}{Theorem}
\newtheorem{mcor}[main]{Corollary}

\newtheorem{theorem}{Theorem}[section]

\newtheorem{lemma}[theorem]{Lemma}
\newtheorem{proposition}[theorem]{Proposition}
\newtheorem{corollary}[theorem]{Corollary}
\theoremstyle{definition}

\newtheorem{remark}[theorem]{Remark}

\begin{document}

\title[Prime II$_1$ factors arising from actions of product groups]
{Prime II$_1$ factors arising from actions of product groups}

\author[D. Drimbe]{Daniel Drimbe}
\address{Department of Mathematics, University of Regina, 3737 Wascana Pkwy, Regina, SK S4S 0A2, Canada.}
\email{daniel.drimbe@uregina.ca}

\begin{abstract} 
We prove that any II$_1$ factor arising from a free ergodic probability measure preserving action $\Gamma\car X$ of a product $\Gamma=\Gamma_1\times\dots\times\Gamma_n$ of icc hyperbolic, free product or wreath product groups is prime, provided $\Gamma_i\car X$ is ergodic, for any $1\leq i\leq n.$ We also completely classify all the tensor product decompositions of a II$_1$ factor associated to a free ergodic probability measure preserving action of a product of icc, hyperbolic, property (T) groups. As a consequence, we derive a unique prime factorization result for such II$_1$ factors. Finally, we obtain a unique prime factorization theorem for a large class of II$_1$ factors  which have property Gamma.

\end{abstract}

\maketitle

\section{Introduction}

\subsection{Background}
In their pioneering work \cite{MvN36,MvN43}, Murray and von Neumann found a natural way to associate a II$_1$ factor, denoted $L(\Gamma)$, to every countable infinite conjugacy class group $\Gamma$ and a II$_1$ factor, denoted $L^\infty(X)\rtimes\Gamma$, to any free ergodic probability measure preserving action $\Gamma\car (X,\mu).$ The classification of these group and group measure space von Neumann algebras is in general a very difficult problem. Nevertheless, a plethora of remarkable results have been obtained in the last 15 years due to S. Popa's influential deformation/rigidity theory, see the surveys \cite{Po07,Va10a, Io12a,Io17}.

A central theme is the study of tensor product decompositions. A II$_1$ factor $M$ is called {\it prime} if it cannot be decomposed as a tensor product of II$_1$ factors. The uncovering of primeness results has been initially explored in the group von Neumann algebra setting. In \cite{Po83}, S. Popa has discovered the first examples of prime II$_1$ factors by showing that the von Neumann algebra of any free group on uncountable many generators is prime. Using D. Voiculescu's free probability theory, L. Ge provided the first examples of separable prime II$_1$ factors by proving that the free group factors $L(\mathbb F_n),$ $2\leq n\leq \infty$, are also prime \cite{Ge96}. By providing new methods in the C$^*$-algebraic setting, N. Ozawa proved that any infinite conjugacy class (icc) hyperbolic group $\Gamma$ gives rise to a {\it solid} II$_1$ factor $L(\Gamma)$, meaning that the relative commutant of any diffuse subalgebra of $L(\Gamma)$ is amenable \cite{Oz03}; in particular it follows that $L(\Gamma)$ is prime. In \cite{Pe06}, by developing an innovative technique based on closable derivations, J. Peterson showed primeness of $L(\Gamma)$, for any icc non-amenable group $\Gamma$ which has positive first Betti number. S. Popa then used his deformation/rigidity theory and gave an alternative proof of solidity of $L(\mathbb F_n)$ \cite{Po06b}. 
The intense research activity over the last decade has resulted in many other primeness results, see \cite{Oz04, Po06a, CI08, CH08, Va10b, Bo12, HV12, DI12, CKP14, Ho15}.

In all these results some negative curvature condition on $\Gamma$ is needed, in the form of a geometric assumption (e.g. $\Gamma$ is a hyperbolic group), or a cohomological assumption (e.g. the existence of a certain unbounded quasi-cocycle). Any of these two conditions can be seen as a ``rank one'' property.  
Concerning the primeness problem in the framework of group measure space von Neumann algebras, the techniques presented in the aforementioned papers can be used to show that any free ergodic probability measure preserving (pmp) action of such groups gives rise to a prime II$_1$ factor. 
Specifically, N. Ozawa showed that $L^\infty(X)\rtimes\Gamma$ is prime whenever $\Gamma\car (X,\mu)$ is a free ergodic pmp action of a non-elementary
hyperbolic group \cite{Oz04} (see also \cite{CS11}). By obtaining new Bass-Serre type rigidity results for II$_1$ factors, I. Chifan and C. Houdayer showed that the II$_1$ factor associated to any free ergodic pmp action of a free product group is prime \cite{CH08}. Then by developing methods from \cite{Si10,Va10b}, D. Hoff proved that $L^\infty(X)\rtimes\Gamma$ is prime whenever $\Gamma\car (X,\mu)$ is a free ergodic pmp action of a group which has positive first Betti number \cite{Ho15}.

\subsection{Statement of the main results}
The first primeness results for group von Neumann algebras arising from icc irreducible lattices in higher rank semisimple Lie groups were obtained only recently in our joint work with D. Hoff and A. Ioana \cite{DHI16} (see also \cite{CdSS17, dSP18}).
Recall that a lattice $\Gamma$ in a product $G=G_1\times\dots\times G_n$ of locally compact second countable groups is called {irreducible} if the action of $G$ on the homogeneous space $G/\Gamma$ is {\it irreducible}, meaning
$G_i\car G/\Gamma$ is ergodic for any $1\leq i\leq n.$ More generally, a pmp action $G\car (X,\mu)$ is called {\it irreducible} if $G_i\car (X,\mu)$ is ergodic for any $1\leq i\leq n.$

Despite all these advancements, the primeness problem for II$_1$ factors arising from arbitrary free ergodic pmp actions of groups of ``higher rank type'' is largely open. Our results aim in this direction by finding a large class of product groups for which all their irreducible actions give rise to prime II$_1$ factors, see Corollary \ref{A2}.
These examples follow from our main technical result. Before stating the result, we introduce the following class of groups and explain the terminology that will be used. 

{\bf Class $\mathcal C$.} We say that a countable group $\Gamma$ belongs to the class $\mathcal C$ if one of the following conditions is satisfied:
\begin{enumerate}
\item $\Gamma$ is an icc, weakly amenable, bi-exact group (see \cite{PV11} for terminology), or
\item $\Gamma=\Sigma_1*\Sigma_2$ is a free product of arbitrary groups such that $|\Sigma_1|\ge 2$ and $|\Sigma_2|\ge 3$, or
\item $\Gamma=\Sigma_0\wr \Gamma_0$ is the wreath product between a non-trivial amenable group $\Sigma_0$ and a non-amenable group $\Gamma_0.$ 
\end{enumerate}

The symbol $\prec$ stands for Popa's intertwining-by-bimodules technique (see Section \ref{cornerr}). We denote by $M^t$ the {\it amplification} of the II$_1$ factor $M$ by $t>0$ and for a pmp action $\Gamma\car (X,\mu)$ we denote by $L^\infty(X)^{\Sigma}$ the subalgebra of elements of $L^\infty(X)$ fixed by a subgroup $\Sigma$ of $\Gamma$ (see Section \ref{term}).

\begin{main}\label{A}
Let $\Gamma=\Gamma_1\times\dots \times\Gamma_n$ be a product of $n\ge 2$ groups that belong to the class $\mathcal C$. Let $\Gamma\overset{}{\car} (X,\mu)$ be a free ergodic pmp action and 
denote $M=L^\infty (X)\rtimes\Gamma.$\\
Suppose that $M= P_1\bar\otimes P_2 $, for some II$_1$ factors $P_1$ and $P_2.$

Then there exists a partition $T_1\sqcup T_2=\{1,\dots ,n \}$ such that $L^\infty(X)\prec_M L^\infty(X)^{\Gamma_{T_1}}\vee L^\infty(X)^{\Gamma_{T_2}}$, where $
\Gamma_{T_i}:=\times_{j\in T_i}\Gamma_j$, for any $i\in\{1,2\}$.

Moreover, if in addition the groups $\Gamma_i$'s have Kazhdan's property (T), then there exist a decomposition $M=P_1^t\bar\otimes P_2^{1/t}$, for some $t> 0$, and a unitary $u\in M$ such that
$$
P_1^t=u (L^\infty(X)^{\Gamma_{T_2}}\rtimes\Gamma_{T_1})u^* \text{ and } P_2^{1/t}=u(L^\infty(X)^{\Gamma_{T_1}}\rtimes\Gamma_{T_2})u^*.
$$

In particular, there exists a pmp action $\Gamma_{T_i}\car (X_i, \mu_i)$ for any $i\in\{1,2\}$ such that
$\Gamma\car X$   is isomorphic to the product action  $ \Gamma_{T_1}\times \Gamma_{T_2}\car X_1\times X_2 .
$

\end{main}

The moreover part applies if the groups $\Gamma_i$'s are icc, property (T), weakly amenable, bi-exact. The following classes of groups satisfy these conditions.

\begin{enumerate}
\item uniform lattices in $Sp(k,1)$ with $k\ge 2$ or any icc group in their measure equivalence class,
\item Gromov's random groups with density satisfying $3^{-1}<d<2^{-1}.$
\end{enumerate}

Note that the moreover part of Theorem \ref{A} provides the first class of product groups for which the primeness problem for II$_1$ factors arising from their actions is completely settled.

\begin{mcor}\label{A2}
Let $\Gamma=\Gamma_1\times\dots \times\Gamma_n$ be a product of $n\ge 2$ groups\footnote{The case $n=1$ already follows from \cite{Oz04}, \cite{CH08} and \cite{CPS11}.} that belong to the class $\mathcal C$. Let $\Gamma\overset{}{\car} (X,\mu)$ be a free ergodic pmp action.

If $\Gamma\car (X,\mu)$ is irreducible, or if the groups $\Gamma_i$'s have Kazhdan's property (T) and the action $\Gamma\car (X,\mu)$ does not admit a direct product decomposition, then $L^\infty(X)\rtimes\Gamma$ is prime.

\end{mcor}

Theorem \ref{A} allows us to prove a unique prime factorization theorem for any II$_1$ factor arising from an arbitrary free ergodic pmp action of a product of groups that belong to the class $\mathcal C$ and have Kazhdan's property (T). More precisely, we have:

\begin{mcor}\label{mcor}\label{C}
Let $\Gamma=\Gamma_1\times\dots \times\Gamma_n$ be a product of $n\ge 2$
groups that belong to the class $\mathcal C$ and have Kazhdan's property (T). Let $\Gamma\overset{}{\car} (X,\mu)$ be a free ergodic pmp action. Denote $M=L^\infty (X)\rtimes\Gamma.$

Then there exists a unique partition $S_1\sqcup \dots \sqcup S_k =\{1,...,n\}$, for some $1\leq k\leq n$, (up to a permutation) and a pmp action $\Gamma_{S_{i}}\car (X_i,\mu_i)$, for any $1\leq i\leq k$, such that:

\begin{enumerate}
\item $\Gamma\car X$ is isomorphic to the product action $\Gamma_{S_1}\times ...\times \Gamma_{S_k}\car X_1\times ...\times X_k$.
\item $M_i:=L^\infty(X_i)\rtimes\Gamma_{S_i}$ is prime for any $1\leq i\leq k$.
\end{enumerate}
Moreover,  the following hold:
\begin{enumerate}
\item If $M=P_1\bar{\otimes}P_2$, for some II$_1$ factors $P_1, P_2$, then there exist a partition $I_1\sqcup I_2=\{1,...,k\}$ and a decomposition $M=P_1^t\bar{\otimes}P_2^{1/t}$, for some $t>0$, such that $P_1^t=\bar{\otimes}_{i\in I_1}M_i$ and $P_2^{1/t}=\bar{\otimes}_{i\in I_2}M_i$, up to unitary conjugacy in $M$.
\item If $M=P_1\bar{\otimes}\dots\bar{\otimes}P_m$, for some $m\geq k$ and II$_1$ factors $P_1,...,P_m$, then $m=k$ and there exists a decomposition $M=P_1^{t_1}\bar{\otimes}...\bar{\otimes}P_k^{t_k}$ for some $t_1,...,t_k>0$ with $t_1t_2\dots t_k=1$ such that after permutation of indices and unitary conjugacy we have $M_i=P_i^{t_i}$, for all $1\leq i\leq k$.

\item In (2), the assumption $m\geq k$ can be omitted if each $P_i$ is assumed to be prime.
\end{enumerate}
\end{mcor}

The first unique prime factorization results for II$_1$ factors were obtained by N. Ozawa and S. Popa in their seminal work \cite{OP03}. Subsequently, several other unique prime factorization results have been obtained in \cite{Pe06, CS11, SW11,Is14, CKP14, HI15, Ho15, Is16, DHI16, De19}. Corollary \ref{C} is the first unique prime factorization result that applies to II$_1$ factors arising from arbitrary free ergodic pmp actions of product groups.

Note that all the known unique prime factorization results in the II$_1$ factor framework are using von Neumann algebras which do not have Murray and von
Neumann's {\it property Gamma} \cite{MvN43}. 
Another novel aspect of this paper is the following unique prime factorization theorem in which the factors possibly have property Gamma. 

\begin{main}\label{UPFgeneral}
For any $1\leq i\leq k$, let $\Gamma_i\car (X_i,\mu_i)$ be a free ergodic pmp action of a group $\Gamma_i$ that belongs to the class $\mathcal C$. For any $1\leq i\leq k$, denote $M_i=L^\infty(X_i)\rtimes\Gamma_i$ and let $M=M_1\bar\otimes \dots \bar\otimes M_k$. Then the following hold:

\begin{enumerate}
\item If $M=P_1\bar{\otimes}P_2$, for some II$_1$ factors $P_1, P_2$, then there exists a partition $I_1\sqcup I_2=\{1,...,k\}$ such that $P_j$ is stably isomorphic to $\bar{\otimes}_{i\in I_j}M_i$, for any $j\in\{1,2\}.$
\item If $M=P_1\bar{\otimes}\dots\bar{\otimes}P_m$, for some $m\geq k$ and II$_1$ factors $P_1,...,P_m$, then $m=k$ and there exists a permutation $\sigma$ of $\{1,...,k\}$ such that $P_i$ is stably isomorphic to $M_{\sigma(i)}$, for any $1\leq i\leq k.$

\item In (2), the assumption $m\geq k$ can be omitted if each $P_i$ is assumed to be prime.
\end{enumerate}
Moreover, if $\Gamma_i\car (X_i,\mu_i)$ is strongly ergodic, for any $1\leq i\leq k$, then the identifications of the von Neumann algebras in (1), (2) and (3) are implemented up to amplification by a unitary from $M$ (as in Corollary \ref{C}).
\end{main}

\begin{remark}
Assume $\Gamma_i\car (X_i,\mu_i)$ is not strongly ergodic, for any $1\leq i\leq k$. Then, the conclusion of Theorem \ref{UPFgeneral} is optimal in the sense that it cannot be improved to deduce that the identifications in (1), (2) and (3) can be implemented up to amplification by a unitary from $M$. 
This follows from \cite[Theorem B]{Ho15}, since all the $M_i$'s have property Gamma.
To ilustrate this, if we assume that $k=2$, \cite[Theorem B]{Ho15} implies that $M$ admits an automorphism $\theta$ such that $\theta(M_i^t)$ is not unitarily conjugate to $M_j$, for any $i,j\in\{1,2\}$ and $t>0.$

\end{remark}


{\bf Comments on the proof of Theorem \ref{A}.} We end the introduction with some informal and brief comments on the proof of Theorem \ref{A}.
For simplicity, assume $\Gamma=\Gamma_1\times...\times \Gamma_n$ is a product of $n\ge 2$ icc, weakly amenable, bi-exact groups.
Let $\Gamma\car (X,\mu)$ be a free ergodic pmp action and denote $M=L^\infty(X)\rtimes\Gamma.$ Assume that we have the tensor product decomposition $M=P_1\bar\otimes P_2$ into II$_1$ factors. We aim to show that $\Gamma\car X$ admits a non-trivial direct product decomposition.

In order to attain this goal we will heavily use S. Popa's deformation/rigidity theory. In the first part of the proof we use S. Popa and S. Vaes' breakthrough work \cite{PV11,PV12} to obtain a partition $T_1\sqcup T_2=\{1,\dots,n\}$ such that 
\begin{equation}\label{z1}
P_1\prec L^\infty(X)\rtimes\Gamma_{T_1}, \text{    and    }  P_2\prec L^\infty(X)\rtimes\Gamma_{T_2},
\end{equation}
where $\Gamma_{T_i}=\times_{j\in T_i}\Gamma_j$, for any $i\in \{1,2\}$. Here, $P\prec Q$ denotes the fact that a corner of $P$ embeds into a corner of $Q$ inside the ambient algebra in the sense of Popa \cite{Po03}. For ease of notation, we will write $P\sim Q$ if $Pp'\prec Q$ and $Qq'\prec P$, for any $p'\in P'\cap M$ and $q'\in Q'\cap M$. 

Since the equality $P_i\vee (P_i'\cap M)= M$ can be seen as a finite index inclusion of von Neumann algebras in the sense of Popa-Pimsner \cite{PP86} for any $i\in\{1,2\}$ (see Section \ref{S: PP}), we can make use of \eqref{z1} and deduce the existence of some abelian von Neumann subalgebras $D_1\subset P_1$ and $D_2\subset P_2$ such that
\begin{equation}\label{z2}
L^\infty(X)\rtimes\Gamma_{T_1}\prec P_1\bar\otimes D_2,    \text{   and    }  L^\infty(X)\rtimes\Gamma_{T_2}\prec D_1\bar\otimes P_2,
\end{equation}
\begin{equation}\label{z3}
D_2\sim {L^\infty(X)}^{\Gamma_{T_1}},   \text{   and   } D_1\sim {L^\infty(X)}^{\Gamma_{T_2}}.
\end{equation}
Here, we denote by ${L^\infty(X)}^{\Gamma_{T_1}}$ and ${L^\infty(X)}^{\Gamma_{T_2}}$ the subalgebras of elements in ${L^\infty(X)}$ fixed by $\Gamma_{T_1}$ and $\Gamma_{T_2},$ respectively.
By combining the intertwining relations \eqref{z2} and \eqref{z3} we show that ${L^\infty(X)}\sim {L^\infty(X)}^{\Gamma_{T_1}}\vee {L^\infty(X)}^{\Gamma_{T_2}}$. 

Finally, if we assume in addition that the groups $\Gamma_i$'s  have property (T), we deduce that we have the identifications
$$
P_1 = {L^\infty(X)}^{\Gamma_{T_2}}\rtimes{\Gamma_{T_1}},
\text{    and   }
P_2 = {L^\infty(X)}^{\Gamma_{T_1}}\rtimes{\Gamma_{T_2}},
$$
up to a unitary conjugacy and amplification.

{\bf Acknowledgment.} I warmly thank Ionut Chifan and Adrian Ioana for many comments and suggestions that helped improve the exposition of the paper.
I am especially grateful to Adrian Ioana for valuable comments on a previous draft which helped increase the generality of the results. I also thank Mart\'in Argerami and Remus Floricel for a useful discussion about these results. Finally, I would like to thank the referee for valuable comments. The author was partially supported by PIMS fellowship.

\section{Preliminaries}

\subsection{Terminology}\label{term} In this paper we consider {\it tracial von Neumann algebras} $(M,\tau)$, i.e. von Neumann algebras $M$ equipped with a faithful normal tracial state $\tau: M\to\mathbb C.$ This induces a norm on $M$ by the formula $\|x\|_2=\tau(x^*x)^{1/2},$ for all $x\in M$. We will always assume that $M$ is a {\it separable} von Neumann algebra, i.e. the $\|\cdot\|_2$-completion of $M$ denoted by $L^2(M)$ is separable as a Hilbert space.
We denote by $\mathcal U(M)$ the {\it unitary group} of $M$ and by $\mathcal Z(M)$ its {\it center}.

All inclusions $P\subset M$ of von Neumann algebras are assumed unital. We denote by $e_P: L^2(M)\to L^2(P)$ the orthogonal projection onto $L^2(P)$, by $E_{P}:M\to P$ the unique $\tau$-preserving {\it conditional expectation} from $M$ onto $P$, by $P'\cap M=\{x\in M|xy=yx, \text{ for all } y\in P\}$ the {\it relative commutant} of $P$ in $M$ and by $\mathcal N_{M}(P)=\{u\in\mathcal U(M)|uPu^*=P\}$ the {\it normalizer} of $P$ in $M$.  We say that $P$ is {\it regular} in $M$ if the von Neumann algebra generated by $\mathcal N_M(P)$ equals $M$.
For two von Neumann subalgebras $P,Q\subset M$, we denote by $P\vee Q$ the von Neumann algebra generated by $P$ and $Q$. {\it Jones' basic construction} of the inclusion $P\subset M$ is defined as the von Neumann subalgebra of $\mathbb B(L^2(M))$ generated by $M$ and $e_P$, and is denoted by $\langle M,e_P \rangle$.

The {\it amplification} of a II$_1$ factor $(M,\tau)$ by a positive number $t$ is defined to be $M^t=p(\mathbb B(\ell^2(\mathbb Z))\bar\otimes M)p$, for a projection $p\in \mathbb B(\ell^2(\mathbb Z))\bar\otimes M$ satisfying $($Tr$\otimes\tau)(p)=t$. Here Tr denotes the usual trace on $\mathbb B(\ell^2(\mathbb Z))$. Since $M$ is a II$_1$ factor, $M^t$ is well defined. Note that if $M=P_1\bar\otimes P_2$, for some II$_1$ factors $P_1$ and $P_2$, then there exists a natural identification $M=P_1^t\bar\otimes M_2^{1/t}$, for every $t>0.$

Let $\Gamma\overset{\sigma}{\car} A$ be a trace preserving action of a countable group $\Gamma$ on a tracial von Neumann algebra $(A,\tau)$. For a subgroup $\Sigma<\Gamma$ we denote by $A^\Sigma=\{a\in A|\sigma_g(a)=a, \text{ for all } g\in\Sigma\}$, the subalgebra of elements of $A$ fixed by $\Sigma.$

Finally, for a product group $\Gamma=\Gamma_1 \times\dots\times \Gamma_n$ and a subset $T\subset \{1,\dots,n\}$, we denote $\Gamma_T=\times _{i\in T}\Gamma_i.$

\subsection {Intertwining-by-bimodules}\label{cornerr} We next recall from \cite [Theorem 2.1 and Corollary 2.3]{Po03} the {\it intertwining-by-bimodules} technique of S. Popa, which gives a powerful criterion for the existence of intertwiners between arbitrary subalgebras of a tracial von Neumann algebra.

\begin {theorem}[\!\!\cite{Po03}]\label{corner} Let $(M,\tau)$ be a tracial von Neumann algebra and let $P\subset pMp, Q\subset qMq$ be von Neumann subalgebras. 
Let $\mathcal G\subset\mathcal U (P)$ be a subgroup such that $\mathcal G''=P$,

Then the following are equivalent:

\begin{itemize}

\item There exist projections $p_0\in P, q_0\in Q$, a $*$-homomorphism $\theta:p_0Pp_0\rightarrow q_0Qq_0$  and a non-zero partial isometry $v\in q_0Mp_0$ such that $\theta(x)v=vx$, for all $x\in p_0Pp_0$.

\item There is no sequence $(u_n)_n\subset\mathcal G$ satisfying $\|E_Q(xu_ny)\|_2\rightarrow 0$, for all $x,y\in M$.

\end{itemize}
\end{theorem}

If one of these equivalent conditions holds true, we write $P\prec_{M}Q$, and say that {\it a corner of $P$ embeds into $Q$ inside $M$.}
If $Pp'\prec_{M}Q$ for any non-zero projection $p'\in P'\cap pMp$, then we write $P\prec^{s}_{M}Q$.\\
Whenever the ambient algebra $(M,\tau)$ is clear from the context, we will write $P\prec Q$ instead of $P\prec_{M}Q$.

The following lemma is a consequence of \cite[Lemma 4.11]{OP07}. For completeness, we provide a short proof.

\begin{lemma}[\!\!\cite{OP07}]\label{L:free}
Let $\Gamma{\car} (Y,\nu)$ be an ergodic pmp action of an icc group and denote $M=L^\infty(Y)\rtimes\Gamma$. If $L^\infty(Y)'\cap M\prec_M L^\infty(Y)$, then $\Gamma\car (Y,\nu)$ is free.
\end{lemma}

{\it Proof.} Let $B=L^\infty(Y)$. The assumption implies that there exist non-zero projections $p\in B'\cap M, q\in B$, a non-zero partial isometry $v\in qMp$ and an injective $*$-homomorphism $\theta: p(B'\cap M)p\to Bq$ such that $\theta(x)v=vx$, for all $x\in p(B'\cap M)p$. 

Note that a standard argument which goes back to \cite{MvN43} shows that $M$ is a factor since $\Gamma\car Y$ is ergodic and $\Gamma$ is icc.
Since $p(B'\cap M)p$ and $Bp$ are abelian, it follows that $p(B'\cap M)p$ is maximal abelian in $pMp$. Using that $M$ is a II$_1$ factor, we obtain that there exists a maximal abelian subalgebra $C\subset M$ such that $Cp=p(B'\cap M)p$ and $p\in C.$ Hence, $C\prec_M B$. We can now apply \cite[Lemma 4.11]{OP07} and obtain the conclusion.
\hfill$\blacksquare$

We continue by observing some elementary facts. The first result is well known and we include a short proof for the reader's convenience.

\begin{lemma}\label{small}
Let $N$ be a von Neumann subalgebra of a tracial von Neumann algebra $(M,\tau)$. Let $P\subset pNp$ and $Q\subset qNq$ be von Neumann subalgebras such that $P\prec_M Q$ and assume that $Q\subset qMq$ is regular.
Then $P\prec_N Q.$
\end{lemma}

{\it Proof.} Assume the contrary, that $P\nprec_N Q$. Thus, there exists a sequence of unitaries $(u_n)_n\subset \mathcal U(P)$ such that $\|E_{Q}(xu_ny)\|_2\to 0$, for any $x,y\in N.$ Thus, $\|E_{Q}(u_ny)\|_2=\|E_{Q}(u_nE_{N}(y))\|_2\to 0$, for any $y\in M.$ Since $Q$ is regular in $qMq$, we obtain that $\|E_{Q}(xu_ny)\|_2\to 0$, for any $x,y\in M$, contradiction.
\hfill$\blacksquare$

\begin{lemma}\label{L: joint}
Let $(M,\tau)$ be a tracial von Neumann algebra and let $Q\subset M$ be a regular von Neumann subalgebra. Let $R_1,R_2\subset M$ be commuting von Neumann subalgebras such that $R_i\prec_M^s Q$ for any $i\in\{1,2\}$. Suppose $R_2$ is abelian. 

Then $R_1\vee R_2\prec^s_M Q.$ 

\end{lemma}

{\it Proof.}
Take a non-zero projection $s\in (R_1\vee R_2)'\cap M$. Since $R_1s\prec_M Q, $ there exist projections $r_1\in R_1$ and $q\in Q, $ a $*$-homomorphism $\theta:r_1R_1r_1s\to qQq$ and a non-zero partial isometry $v\in qMr_1s$, satisfying 
\begin{equation}\label{111}
\theta(x)v=vx, \text{   for all   } x\in r_1R_1r_1s.
\end{equation}

We argue that $(r_1R_1r_1s)\vee (R_2r_1s)\prec_M Q.$ By supposing the contrary, there exist two sequences of unitaries $(u_n)_n\subset \mathcal U(r_1R_1r_1s)$ and $(v_n)_n\subset \mathcal U(R_2)$ such that 
$$
\|E_{Q}(xu_n (v_nr_1s)y)\|_2\to 0, \text{   for all    } x,y\in M.
$$

Note that $u_nr_1s=u_n$, for all $n$.
By taking $x=v$ and using the intertwining relation \eqref{111}, we get that 
$$
\|E_{Q}(\theta(u_n)v  v_n y)\|_2=\|E_{Q}(v v_n y)\|_2\to 0, \text{   for all    } y\in M.
$$
Let $r:=v^*v.$ Since $Q$ is regular, we get that
\begin{equation}\label{zz}
\|E_{Q}(xr v_n y)\|_2\to 0, \text{   for all    } x,y\in M.
\end{equation}
Denote $r':=\vee_{w\in\mathcal U(R_2)}wrw^*\in R_2'\cap M$. Since $R_2$ is abelian, \eqref{zz} implies that $\|E_{Q}(x(wrw^*) v_n y)\|_2\to 0,$ for any $w\in\mathcal U(R_2)$ and $x,y\in M.$
Note that $p_1\vee p_2=s(p_1+p_2)$, for any two projections $p_1$ and $p_2$ in $M$. Here we denote by $s(b)$ the support projection of a positive element $b\in M.$ 
Moreover, by using Borel functional calculus, there exist a sequence $(c_n)_n\subset M$ such that $c_n(p_1+p_2)$ converges to $s(p_1+p_2)$ in the $\|\cdot\|_2$-norm.
Therefore, it follows that $\|E_{Q}(x (w_1rw_1^*+w_2rw_2^*) v_n y)\|_2\to 0$, and hence, $\|E_{Q}(x (w_1rw_1^*\vee w_2rw_2^*) v_n y)\|_2\to 0$, for any $w_1,w_2\in \mathcal U(R_2)$ and $x,y\in M$.

Finally, by induction it follows that 
$$
\|E_{Q}(x(w_1rw_1^*\vee...\vee w_mrw_m^*) v_n y)\|_2\to 0,
$$
for any $w_1,...,w_m\in\mathcal U(R_2)$ and $x,y\in M$. \\
Hence, $\|E_{Q}(xr'v_ny)\|_2\to 0$, for all $x,y\in M$. This shows that $R_2r'\nprec_M Q,$ contradiction. Therefore, $(R_1\vee R_2)s\prec_M Q$, which implies that $R_1\vee R_2\prec_M^s Q.$
\hfill$\blacksquare$

We will need the following result which is an extension of \cite[Lemma 2.8(2)]{DHI16} (see also \cite[Proposition 2.7]{PV11}).

\begin{proposition}\label{L: PV}
Let $(M,\tau )$ be a tracial von Neumann algebra and let $Q_1,Q_2\subset M$ be von Neumann subalgebras which form a commuting square, i.e. $E_{Q_1}\circ E_{Q_2}=E_{Q_2}\circ E_{Q_1}$. Assume that there exist commuting subgroups $\mathcal N_1< \mathcal N_M(Q_1)$ and $\mathcal N_2< \mathcal N_M(Q_2)$ satisfying $(\mathcal N_1\vee\mathcal N_2)''=M$. 
Let $P\subset pMp$ be a von Neumann subalgebra.

If $P\prec_M^s Q_1$ and $P\prec_M^s Q_2$, then $P\prec_M^s Q_1\cap Q_2.$
\end{proposition}

\begin{remark}
The proposition will be applied for the following particular case. Assume $M=P_1\bar\otimes P_2$ for some von Neumann subalgebras $P_1,P_2\subset M$ and let $D_i\subset P_i$ be a subalgebra for any $i\in\{1,2\}.$ By taking  $Q_1=D_1\bar\otimes P_2$ and $Q_2=P_1\bar\otimes D_2$ the assumptions of Proposition \ref{L: PV} are satisfied.
\end{remark}

The proof of Proposition \ref{L: PV} follows directly by using the next lemma and adapting the proof of \cite[Lemma 2.8(2)]{DHI16} (see also \cite[Proposition 2.7]{PV11}). We leave the other details to the reader.

\begin{lemma}
Let $(M,\tau)$ be a tracial von Neumann algebra and let $Q_1,Q_2\subset M$ and $\mathcal N_1,\mathcal N_2\subset M$ be as in Proposition \ref{L: PV}. Denote $Q=Q_1\cap Q_2.$

Then the $M$-$M$-bimodule $L^2(\langle M,e_{Q_1} \rangle)\otimes_M L^2(\langle M,e_{Q_2} \rangle)$ is contained in a multiple of the $M$-$M$-bimodule $L^2(\langle M,e_{Q}\rangle)$.
\end{lemma}

{\it Proof.} The proof follows almost verbatim part of the proof of \cite[Proposition 2.7]{PV11}. However, we provide some details for the reader's convenience. 

Denote by $H$ the $M$-$M$-bimodule $L^2(\langle M,e_{Q_1} \rangle)\otimes_M L^2(\langle M,e_{Q_2} \rangle)$. For $u_1,v_1\in \mathcal N_1$ and $u_2,v_2\in\mathcal N_2$, denote by $H_{u_1,v_1}^{u_2,v_2}$ the closed linear span of $\{xe_{Q_1}u_1u_2\otimes_M v_1v_2e_{Q_2}y|x,y\in M\}.$ Note that the formulas $u E_{Q_i}(\cdot)u^*=E_{Q_i}(u\cdot u^*)$, for any $u\in\mathcal N_M(Q_i)$, combined with the commuting square property imply that the map
$$
x{e_{Q_1}}u_1u_2\otimes_M v_1v_2e_{Q_2}y\to xu_1v_1\otimes_Q u_2v_2y
$$
defines an $M$-$M$-bimodular unitary operator of $H_{u_1,v_1}^{u_2,v_2}$  to $L^2(\langle M,e_{Q}\rangle)$. To show this it suffices to verify that
\begin{equation}\label{formula}
\langle {xe_{Q_1}}u_1u_2\otimes_M v_1v_2e_{Q_2}y, {e_{Q_1}}u_1u_2\otimes_M v_1v_2e_{Q_2} \rangle=\langle x u_1v_1\otimes_Q u_2v_2y,u_1v_1\otimes_Q u_2v_2 \rangle,
\end{equation}
for all $x,y\in M$. Note that the left hand side of \eqref{formula} equals
$$
\begin{array}{rcl}
\langle 
{xe_{Q_1}}u_1u_2\otimes_M v_1v_2e_{Q_2}y, {e_{Q_1}}u_1u_2\otimes_M v_1v_2e_{Q_2} 
\rangle
&=&
 \tau(E_{Q_2}(v_2^*v_1^*  E_M(u_2^*u_1^*E_{Q_1}(x)u_1u_2)  v_1v_2)y)\\
&=& \tau(v_2^* u_2^*E_{Q_2}(v_1^* u_1^*E_{Q_1}(x)u_1 v_1)u_2 v_2y)
\end{array}
$$
Therefore, the right hand side or \eqref{formula} equals its left hand side since
$$
\langle x u_1v_1\otimes_Q u_2v_2y,u_1v_1\otimes_Q u_2v_2 \rangle 
= \tau(v_2^*u_2^*E_{Q_2}(E_{Q_1}(v_1^*u_1^*xu_1v_1))u_2v_2y).
$$

Remark that the regularity assumption on the $Q_i$'s implies that 
the closed linear span of $\{u_1u_2|u_1\in\mathcal N_1, u_2\in\mathcal N_2\}$ equals to $L^2(M)$.
Therefore, the closed linear span of 
$\{H_{u_1,v_1}^{u_2,v_2}|$ $u_1,v_1\in \mathcal N_1,u_2,v_2\in \mathcal N_2\}
$ 
equals to $H$. This shows that $H$ is contained in a multiple of $L^2(\langle M,e_Q\rangle)$. 
\hfill$\blacksquare$

\subsection{Finite index inclusions of von Neumann algebras.}\label{S: PP}

For an inclusion $P\subset M$ of II$_1$ factors the {\it Jones index} is the dimension of $L^2(M)$ as a left $P$-module \cite{Jo81}. In \cite{PP86}, M. Pimsner and S. Popa defined a probabilistic notion of index for an inclusion $P\subset M$ of arbitrary von Neumann algebras with conditional expectation, which in the case of inclusions of II$_1$ factors coincides with Jones' index.
Following \cite{PP86}, we say that the inclusion $P\subset M$ of tracial von Neumann algebras has {\it probabilistic index} $[M:P]=\lambda^{-1}$, where
$$
\lambda= \text{inf}\{\|E_P(x)\|^2_2{\|x\|_2^{-2}}|x\in M_+, x\neq 0\}.
$$
Here we use the convention that $\frac{1}{0}=\infty.$

\begin{lemma}[\!\!{\cite[Lemma 2.3]{PP86}}]\label{PP}
Let $P\subset M$ be an inclusion of tracial von Neumann algebras such that $[M:P]<\infty$. Then the following hold:
\begin{enumerate}
\item If $p\in P$ is a projection, then $[pMp: pPp]<\infty.$
 
\item $M\prec^s_M P.$
\end{enumerate}
\end{lemma}

For a proof, see \cite[Lemma 2.4]{CIK13}. We will need the following well known lemma and we include its proof for completeness (see also \cite[Lemma 3.9]{Va08}).

\begin{lemma}\label{L: center}\label{L: fi}
Let $(M,\tau)$ be a tracial von Neumann algebra and let $R\subset N\subset pMp$ be von Neumann subalgebras such that $[N:R]<\infty$.
Then the following hold:

\begin{enumerate}
\item If $R'\cap N\subset R$, then there exists a non-zero projection $z\in\mathcal Z(R)$ such that $\mathcal Z(N)z=\mathcal Z(R)z$.
\item Assume $R'\cap N\subset R$ or $\mathcal Z(R)$ is completely atomic. If $Q\subset qMq$ is a von Neumann subalgebra such that $R\prec_{M} Q$, then $N\prec_{M} Q.$
\end{enumerate}
\end{lemma}

{\it Proof.} 
(1) Applying Lemma \ref{PP}(2), we get that $N\prec_N R$. By passing to relative commutants, we can use \cite[Lemma 3.5]{Va08} and deduce that $\mathcal Z(R)\prec_N \mathcal Z(N).$ Hence, we obtain that there exist projections $r\in \mathcal Z(R), n\in\mathcal Z(N)$, a non-zero partial isometry $v\in nNr$ and a $*$-homomorphism $\theta : \mathcal Z (R)r\to \mathcal Z(N)n$ such that $v\theta(x)=\theta(x)v=vx$, for all $x\in \mathcal Z (R)r.$ 
By noticing that $\mathcal Z(N)\subset \mathcal Z(R)$, we obtain that $E_{\mathcal Z(R)}(v^*v)\theta(x)=E_{\mathcal Z(R)}(v^*v) x$, for all $x\in \mathcal Z(R)r.$ If we denote by $p_0$ the support projection of $E_{\mathcal Z(R)}(v^*v)$, we get that $p_0\theta(x)=p_0x$, for all $x\in \mathcal Z(R)r.$ Therefore, $\mathcal Z(R)\prec_{\mathcal Z(R)} \mathcal Z(N)$, which clearly implies the conclusion.

(2) Assume first that $R'\cap N\subset R$.
Since $R\prec_M Q$, there exist projections $r\in R$, $q_0\in Q$, a non-zero partial isometry $v\in q_0Mr$ and a $*$-homomorphism $\theta: rRr\to q_0Qq_0$ such that
$$
\theta(x)v=vx, \text{   for all   } x\in rRr.
$$ 
Note that $v^*v\in (R'\cap pMp)r$ and denote by $r_1$ the support projection of $E_{N}(v^*v)$. Notice that $r_1\in (R'\cap N)r.$
Since $R'\cap N\subset R$, then $r_1\in rRr$ and therefore, Lemma \ref{PP} implies that $r_1Nr_1\prec_{r_1Nr_1} r_1Rr_1.$ 

Thus, there exist  projections $n\in r_1Nr_1$, $r_0\in r_1Rr_1$, a non-zero partial isometry $w\in r_0Nn$ and a $*$-homomorphism $\psi: nNn\to r_0Rr_0$ such that 
$$
\psi(x)w=wx, \text{   for all   } x\in nNn.
$$
Moreover, by restricting $ww^*$ if necessary we can assume without loss of generality that the support projection of $E_{r_1Rr_1}(ww^*)$ equals $r_0.$
Note that $\theta(\psi(\cdot)): nNn\to q_0Qq_0$ is a $*$-homomorphism which satisfies 
\begin{equation}\label{c}
\theta(\psi(x))vw=vwx,  \text{   for all   } x\in nNn.
\end{equation}
If $vw=0$, then $E_N(v^*v)ww^*=0$, which implies $r_1ww^*=0.$ Hence, $r_1E_{r_1Rr_1}(ww^*)=0$, showing that $r_0=0$, contradiction. This proves that $vw\neq 0.$ By replacing $vw$ by the partial isometry from its polar decomposition, the intertwining relation \eqref{c} still holds. This shows that $N\prec_M Q.$

Assume now that $\mathcal Z(R)$ is completely atomic. Note that \cite[Lemma 2.11 and Lemma 2.4(2)]{Dr19} give that $N\prec_M Q.$
\hfill$\blacksquare$

\subsection{Relative amenability}
A tracial von Neumann algebra $(M,\tau)$ is {\it amenable} if there exists a positive linear functional $\Phi:\mathbb B(L^2(M))\to\mathbb C$ such that $\Phi_{|M}=\tau$ and $\Phi$ is $M$-{\it central}, meaning $\Phi(xT)=\Phi(Tx),$ for all $x\in M$ and $T\in \mathbb B(L^2(M))$. The famous theorem of A. Connes asserts that a von Neumann algebra $M$ is amenable if and only if it is approximately finite dimensional \cite{Co76}. 

N. Ozawa and S. Popa have considered a very useful relative version of this notion\cite{OP07}. Let $(M,\tau)$ be a tracial von Neumann algebra. Let $p\in M$ be a projection and $P\subset pMp,Q\subset M$ be von Neumann subalgebras. Following \cite[Definition 2.2]{OP07}, we say that $P\subset pMp$ is {\it amenable relative to $Q$ inside $M$} if there exists a positive linear functional $\Phi:p\langle M,e_Q\rangle p\to\mathbb C$ such that $\Phi_{|pMp}=\tau$ and $\Phi$ is $P$-central. 
 Note that $P$ is amenable relative to $\mathbb C$ inside $M$ if and only if $P$ is amenable.

The following lemma is well known and it goes back to \cite[Lemma 10.2]{IPV10}, but we include a proof for completeness.
The arguments are essentially contained in the proof of \cite[Proposition 3.2]{PV12}.

\begin{lemma}\label{ipv}
Let $\Gamma\car (B,\tau)$ be a trace preserving action and denote $M=B\rtimes\Gamma.$ Define the $*$-homomorphism $\Delta:M\to M\bar\otimes L(\Gamma)$ by letting $\Delta(bu_g)=bu_g\otimes u_g,$ for all $b\in B$ and $g\in\Gamma.$\\
Let $P\subset pMp$ be a von Neumann subalgebra such that there exists $p_1\in \Delta(P)'\cap \Delta(p)(M\bar\otimes M)\Delta(p)$ with the property that $\Delta(P)p_1$ is amenable relative to $M\otimes 1$.

Then there exists a non-zero projection $p_0\in P'\cap pMp$ such that $Pp_0$ is amenable relative to $B$ inside $M$.
\end{lemma}

{\it Proof.} Define $\mathcal M=M\bar\otimes M.$ The assumption implies the existence of a positive linear functional $\Phi:p_1\langle \mathcal M, e_{M\otimes 1} \rangle p_1\to \mathbb C$ such that the restriction of $\Phi$ to $p_1\mathcal Mp_1$ equals the trace on $p_1\mathcal Mp_1$ and $\Phi$ is $\Delta(P)p_1$-central. Since $E_{M\bar\otimes 1}\circ\Delta=\Delta\circ E_B$, note that we can define the injective $*$-homomorphism $\Delta_1: \langle M,e_B \rangle\to \langle \mathcal M, e_{M\otimes 1} \rangle$ by letting $\Delta_{1}(e_B)=e_{M\otimes 1}$ and $\Delta_{1}(x)=\Delta(x)$, for all $x\in M$.

Define the positive linear functional $\Psi: p\langle M,e_B\rangle p\to\mathbb C$ by $\Psi(x)=\Phi(p_1\Delta_1(x)p_1)$, for all $x\in p\langle M,e_B\rangle p.$ Note that $\Psi$ is $P$-central and its restriction to $pMp$ is normal. Therefore, \cite[Lemma 2.9]{BV12} implies that there exists a non-zero projection $p_0\in P'\cap pMp$ such that $Pp_0$ is amenable relative to $B$.
\hfill$\blacksquare$

\subsection{Relatively strongly solid groups}

Following \cite[Definition 2.7]{CIK13}, a countable group $\Gamma$ is said to be {\it relatively strongly solid} and write $\Gamma\in \mathcal C _{rss}$ if for any trace preserving action $\Gamma\car B$ the following holds: if $M=B\rtimes\Gamma$ and $A\subset pMp$ is a von Neumann algebra which is amenable relative to $B$, then either $A\prec_M B$ or the normalizer $\mathcal N_{pMp}(A)''$ is amenable relative to $B$.
In their breakthrough work \cite{PV11,PV12}, S. Popa and S. Vaes proved that all non-elementary hyperbolic groups belong to $\mathcal C_{rss}.$ More generally, [PV12, Theorem 1.4] shows that all weakly amenable, bi-exact groups are relatively strongly solid.

A remarkable subsequent development has been made by A. Ioana \cite{Io12} (see also \cite{Va13}) in the context of amalgamated free products by classifying all subalgebras of $M=M_1*_B M_2$ that are amenable relative to $B$ and that satisfy a certain spectral gap condition.

We will make use of the following consequence for groups that belong to $\mathcal C_{\text{rss}}$ (see \cite[Lemma 5.2]{KV15}).

\begin{lemma}[\!\!{\cite{KV15}}]\label{L: rss}
Let $\Gamma\car B$ be a trace preserving action of a group $\Gamma\in \mathcal C_{\text{rss}}$. Denote $M=B\rtimes\Gamma$.
Let $P_1, P_2\subset pM$p be commuting von Neumann subalgebras.

Then either $P_1\prec_{M}B$ or $P_2$ is amenable relative to $B$ inside $M$.
\end{lemma}

For free product groups we will use the following consequence \cite[Theorem 3.1]{CdSS17} of \cite[Theorem A]{Va13}:

\begin{lemma}[\!\!{\cite{CdSS17}}]\label{L: amalgam}
Let $\Gamma\car B$ be a trace preserving action, where $\Gamma=\Sigma_1*\Sigma_2$ and $|\Sigma_1|\ge 2$ and $|\Sigma_2|\ge 3$.  
Denote $M=B\rtimes\Gamma$ and assume $P_1,P_2\subset pMp$ are two commuting diffuse subalgebras such that $P_1\vee P_2\subset pMp$ has finite index. 

Then there exists an $i\in\{1,2\}$ such that $P_i\prec_M B$.
 
\end{lemma}

\subsection{Wreath product groups}

The next lemma gives a dichotomy result for commuting subalgebras of von Neumann algebras arising from trace preserving actions of wreath product groups. The arguments rely heavily on \cite{IPV10}.

\begin{lemma}\label{wreath}
Let $\Gamma=\Sigma_0\wr\Gamma_0$ be the wreath product between a non-trivial amenable group $\Sigma_0$ and an infinite group $\Gamma_0$. Let $\Gamma\car B$ be a trace preserving action and define $M=B\rtimes\Gamma.$ Let $P_1,P_2\subset pMp$ be two commuting subalgebras such that $P_1\vee P_2\subset pMp$ has finite index.

Then there exists a non-zero projection $p_0\in P_1'\cap pMp$ such that $P_1p_0$ is amenable relative to $B$ or $P_2\prec_M B.$

\end{lemma}

{\it Proof.} Let $\Delta: M\to M\bar\otimes L(\Gamma)$ be the $*$-homomorphism defined by $\Delta(bu_g)=bu_g\otimes u_g$, for all $b\in B$ and $g\in\Gamma.$
By applying \cite[Corollary 4.3]{IPV10}, one of the following possibilities occurs: (1) there exists a non-zero projection $p_1\in \Delta(P_1)'\cap \Delta(p)(M\bar\otimes M)\Delta(p)$ such that $\Delta(P_1)p_1$ is amenable relative to $M\otimes 1$, or (2) $\Delta(P_2)\prec M\otimes 1$, or (3) $\Delta(P_1\vee (P_1'\cap pMp))\prec M\bar\otimes L(\Sigma_0^{(\Gamma_0)})$, or (4) $\Delta(P_1\vee (P_1'\cap pMp))\prec M\bar\otimes L(\Gamma_0).$

If (1) holds, by Lemma \ref{ipv} there exists a non-zero projection $p_0\in P_1'\cap pMp$ such that $P_1p_0$ is amenable relative to $B$. If (2) holds, \cite[Lemma 2.9(1)]{Io10} implies that $P_2\prec_M B.$

We end the proof by showing that (3) and (4) cannot hold. Indeed, if (3) or (4) is true, then Lemma \ref{L: fi} combined with \cite[Lemma 9.2(4)]{Io10} imply that $L(\Gamma)\prec_M L(\Sigma_0^{(\Gamma_0)})$ or $L(\Gamma)\prec_M L(\Gamma_0).$ Both are in contradiction with the fact that $\Sigma_0$ and $\Gamma_0$ are infinite groups.
\hfill$\blacksquare$

Combining the previous three lemmas with \cite[Lemma 2.6(2)]{DHI16}, we obtain the following corollary.

\begin{corollary}\label{all}
Let $\Gamma$ be a group that belongs to the class $\mathcal C$. Let $\Gamma\car B$ be a trace preserving action and define $M=B\rtimes\Gamma.$ Let $P_1,P_2\subset pMp$ be two commuting diffuse subalgebras such that $P_1\vee P_2\subset pMp$ has finite index.

Then there exists a non-zero projection $p_0\in \mathcal N_{pMp}(P_1)'\cap pMp$ such that $P_1p_0$ is amenable relative to $B$ or $P_2\prec_M B.$

\end{corollary}

\section{From tensor decompositions of II$_1$ factors to decompositions of actions}
The goal of this section is to prove Theorem \ref{AA} which is our main ingredient of the proof of Theorem \ref{A}.
The moreover part will provide a von Neumann algebraic criterion for pmp actions of product groups to admit a direct product decomposition. First, we need the following result.

\begin{theorem}\label{Th:split}
Let $\Gamma=\Gamma_1\times\Gamma_2$ be a product of countable icc groups and let $\Gamma\overset{}{\car} (X,\mu)$ be a free ergodic pmp action. Denote $M=L^\infty(X)\rtimes\Gamma.$ 

Suppose that $M= P_1\bar\otimes P_2 $ for some II$_1$ factors $P_1$ and $P_2$ such that $P_i\prec_ M^s L^\infty(X)\rtimes\Gamma_i$, for all $i\in\{1,2\}.$

The following hold:

\begin{enumerate}

\item If $L^\infty(X)\prec_{M}L^\infty(X)^{\Gamma_1}\vee L^\infty(X)^{\Gamma_2}$ and $\Gamma_1$  has property (T), 
then there exist a decomposition $M=P_1^t\bar\otimes P_2^{1/t}$, for some $t> 0$, and a unitary $u\in M$ such that
$
P_1^t=u (L^\infty(X)^{\Gamma_{2}}\rtimes\Gamma_{1})u^* \text{ and } P_2^{1/t}=u(L^\infty(X)^{\Gamma_{1}}\rtimes\Gamma_{2})u^*.
$

\item If $L^\infty(X)=L^\infty(X)^{\Gamma_1}\vee L^\infty(X)^{\Gamma_2}$, then $P_1$ is stably isomorphic to $L^\infty(X)^{\Gamma_2}\rtimes\Gamma_1$ and $P_2$ is stably isomorphic to $L^\infty(X)^{\Gamma_1}\rtimes\Gamma_2.$

\end{enumerate}

\end{theorem}

{\it Proof.} Let $\{u_g\}_{g\in\Gamma}$ be the canonical unitaries in $M$ implementing the action $\Gamma\car (X,\mu)$. Denote $A=L^\infty(X)$, $A_1=L^\infty(X)^{\Gamma_2}, A_2=L^\infty(X)^{\Gamma_1}$, $B=A_1\vee A_2$ and $M_i=A_i\rtimes\Gamma_i$, for any $i\in\{1,2\}$. Since $B$ is a $\Gamma$-invariant subalgebra of $A$, we consider $\Gamma\overset {\rho}{\car} B$ the natural action of $\Gamma$ on $B$. Note that \cite[Lemma 3.5]{Va08} shows that $B'\cap (B\rtimes\Gamma)\prec_M A$, which implies by \cite[Lemma 3.7]{Va08} that $B'\cap (B\rtimes\Gamma)\prec_M B.$ Since $B\subset M$ is regular, Lemma \ref{small} shows that $B'\cap (B\rtimes\Gamma)\prec_{B\rtimes\Gamma} B.$ An application of Lemma \ref{L:free} gives us that $\Gamma\car B$ is free. 

Since $\Gamma\overset{\rho}{\car} B$ is isomorphic to the product action $\Gamma_1\times\Gamma_2 \car A_1\bar\otimes A_2$, it follows that $\Gamma_i\car A_i$ is free for any $i\in\{1,2\}$. Hence, 
\begin{equation}\label{fr}
A_2'\cap M=A\rtimes\Gamma_1.
\end{equation}
Indeed, take $x=\sum_{g\in\Gamma}x_gu_g\in A_2'\cap M.$ It follows that $x_gb=x_g\rho_g(b)$, for all $b\in A_2.$ Therefore, $E_{A_2}(x_g^*x_g)b=E_{A_2}(x_g^*x_g)\rho_g(b)$, for all $b\in A_2$. Since $\Gamma_2$ acts freely on $A_2$, we get that $E_{A_2}(x_g^*x_g)=0$, for all $g\notin \Gamma_1$, which implies that $x_g=0$.

Note that $P_1\prec_M A\rtimes\Gamma_{1}$ implies $A_2\prec_M P_2$ by applying \cite[Lemma 3.5]{Va08}.
Since $A_2\subset M$ is regular, we can apply \cite[Corollary 1.3]{Io06} (see also \cite[Proposition 12]{OP07}) and obtain that there exist a unitary $u_0\in\mathcal U(M)$ and a decomposition $M=P_1^{1/t_{0}}\bar\otimes P_2^{t_0}$, for some $t_0>0$ such that $u_0A_2u_0^*\subset P_2^{t_0}.$ Denote $Q_1=u_0^*P_1^{1/t_0}u_0$ and $Q_2=u_0^*P_2^{t_0}u_0.$
Since $A_2\subset Q_2$ we can apply Ge's tensor splitting theorem \cite[Theorem A]{Ge95} and derive that  
\begin{equation}\label{a1}
A\rtimes\Gamma_1=Q_1\bar\otimes (A_2'\cap Q_2).
\end{equation}

Since $A_2'\cap Q_2\subset A\rtimes\Gamma_1$ and $Q_2\prec_M A\rtimes\Gamma_2$, it follows that $A_2'\cap Q_2\prec_M A.$
Thus, there exists a non-zero projection $a\in A_2'\cap Q_2$ such that $a(A_2'\cap Q_2)a$ is abelian. By cutting the von Neumann algebras of the equality \eqref{a1} by the projection $a$ and by passing to the center, we obtain that $A_2a=a(A_2'\cap Q_2)a$ and hence 
\begin{equation}\label{aa1}
a(A\rtimes\Gamma_1)a=Q_1a\bar\otimes A_2a. 
\end{equation}

(1) Now we can prove the first conclusion of the theorem. We first show that $L(\Gamma_1)\prec_M P_1.$ To see this, note that since $L(\Gamma_1)$ is a II$_1$ factor, we obtain by \cite[Lemma 4.5]{CdSS17} that there exists a unitary $u\in A\rtimes\Gamma_1$ such that $b:=uau^*\in L(\Gamma_1).$ Hence, relation \eqref{aa1} gives that $b(A\rtimes\Gamma_1)b=uQ_1u^*b\bar\otimes  uA_2u^*b.$ 
 Note that $bL(\Gamma_1)b$ is a von Neumann algebra with property (T) by \cite[Proposition 4.7(2)]{Po01} and $ u(A_2'\cap Q_2)u^*b$ is an abelian von Neumann algebra. Therefore, we obtain that $bL(\Gamma_1)b\prec_{b(A\rtimes\Gamma_1)b} uQ_1u^*b$ (see e.g. \cite[Lemma 1]{HPV10}). This implies that $L(\Gamma_1)\prec_M P_1.$

By passing to relative commutants and by applying twice\cite[Lemma 3.5]{Va08}, we obtain that $P_2\prec_M M_2$ and hence $M_1\prec_M P_1.$ 
Since $M_1$ and $M_1'\cap M=M_2$ are factors, we can apply \cite[Proposition 12]{OP03} and obtain that there exist a unitary $v\in M$ and a decomposition $M=P_1^s\bar\otimes P_2^{1/s}$ for a positive $s>0$ such that $vM_1v^*\subset P_1^s.$ Therefore, $P_2^{1/s}\subset vM_2v^*.$ By applying \cite[Theorem A]{Ge95}, we obtain that there exists a factor $D\subset P_1^s$ such that $D\bar\otimes P_2^{1/s}=vM_2v^*.$ It is easy to see that $D$ is not diffuse, which implies that $D=\mathbb M_k(\mathbb C) $ for some integer $k\ge 1.$
By denoting $s_0=s/k$, we deduce that $P_2^{1/{s_0}}=vM_2v^*$ and $P_1^{s_0}=vM_1v^*.$ 

(2) Note that the assumption implies $A=A_1\bar\otimes A_2$ and the relation \eqref{aa1} shows that  $a(M_1\bar\otimes A_2)a=Q_1a\bar\otimes A_2a$. Therefore, by disintegrating in the above equality over the center $A_2a$ and by using \cite[Theorem IV.8.23]{Ta01}, we deduce that $M_1$ is stably isomorphic to $Q_1$, hence to $P_1$. In a similar way, we obtain that $M_2$ is stably isomorphic to $P_2$.
\hfill$\blacksquare$

\begin{theorem}\label{AA}
Let $\Gamma=\Gamma_1\times\Gamma_2$ be a product of countable icc groups and let $\Gamma\overset{\sigma}{\car} (X,\mu)$ be a free ergodic pmp action.
Suppose that $M=L^\infty(X)\rtimes\Gamma= P_1\bar\otimes P_2 $ for some II$_1$ factors $P_1$ and $P_2$ such that $P_i\prec_ M^s L^\infty(X)\rtimes\Gamma_i$, for all $i\in\{1,2\}.$

Then $L^\infty(X)\prec_M L^\infty(X)^{\Gamma_1}\vee L^\infty(X)^{\Gamma_2}.$

Moreover, assume that $\Gamma_1$ has property (T). Then there exist a decomposition $M=P_1^t\bar\otimes P_2^{1/t}$, for some $t> 0$, and a unitary $u\in M$ such that
$$
P_1^t=u (L^\infty(X)^{\Gamma_{2}}\rtimes\Gamma_{1})u^* \text{ and } P_2^{1/t}=u(L^\infty(X)^{\Gamma_{1}}\rtimes\Gamma_{2})u^*.
$$

In particular, there exists a pmp action $\Gamma_{i}\car (X_i, \mu_i)$ for any $i\in\{1,2\}$ such that the actions
$\Gamma\car X $ and $ \Gamma_{1}\times \Gamma_{2}\car X_1\times X_2 
$ are isomorphic.

\end{theorem}

{\it Proof.} Let $\{u_g\}_{g\in\Gamma}$ be the canonical unitaries in $M$ implementing the action $\Gamma\car (X,\mu)$ and
denote $A=L^\infty(X)$.
For any $i\in\{1,2\}$, \cite[Proposition 2.4]{CKP14} implies that there exist non-zero projections $p_i\in P_i$, $q_i\in A\rtimes\Gamma_{i}$, a subalgebra $Q_i\subset q_i(A\rtimes\Gamma_{i})q_i$, a partial isometry $v_i\in q_iMp_i$ and a $*$-isomorphism $\theta_i: p_iP_ip_i\to Q_i$ such that:
\begin{equation}\label{s1}
Q_i\vee (Q_i'\cap q_i(A\rtimes\Gamma_{i})q_i)\subset q_i(A\rtimes\Gamma_{i})q_i \text{  has finite index, }
\end{equation}
\begin{equation}\label{s2}
\theta_i(x)v_i=v_ix, \text{   for all   } x\in p_iP_ip_i, \text{  and }
\end{equation}
\begin{equation}\label{s12}
E_{A\rtimes\Gamma_{i}}(v_iv_i^*)q_i\ge \lambda_i q_i,  \text{   for some positive number } \lambda_i.
\end{equation}

The rest of the proof is divided between four claims.

{\bf Claim 1.} We can assume, in addition to \eqref{s1}-\eqref{s12}, that $Q_i$ also satisfies that $Q_i'\cap q_i(A\rtimes\Gamma_{i})q_i=A^{\Gamma_{i}}q_i$, for any $i\in\{1,2\}.$

{\it Proof.} For simplicity we prove the claim only for $i=1$. Denote $R=Q_1'\cap q_1(A\rtimes\Gamma_{1})q_1$. First, note that $R\prec_M A.$ Indeed, on one hand since $P_1\prec_M Q_1$, by considering the relative commutants, \cite[Lemma 3.5]{Va08} implies that $R\prec_M P_2.$ Applying \cite[Lemma 3.7]{Va08} we get $R\prec_M A\rtimes{\Gamma_{2}}$. On the other hand, $R\subset q_1(A\rtimes\Gamma_{1})q_1$. Hence, one can check that we actually have $R\prec_M A.$

Therefore, there exists a non-zero projection $r\in R$ such that $rRr$ is abelian. Applying Lemma \ref{PP}(1), we have that $Q_1r\vee rRr\subset r(A\rtimes\Gamma_{1})r$ has finite index. Since $Q_1$ is a factor, Lemma \ref{L: center}(1) shows that by replacing $r$ by a smaller projection in $R$, we can assume that 
$$rRr=\mathcal Z(rRr)=\mathcal Z(r(A\rtimes{\Gamma_{1}})r)=A^{\Gamma_{1}}r.
$$
Lemma \ref{PP}(1) guarantees that $Q_1r\vee rRr\subset r(A\rtimes{\Gamma_{1}})r$ still has finite index.

Since $E_{A\rtimes\Gamma_{1}}(v_	1v_1^*)q_1\ge \lambda q_1$ and $r\leq q_1$, we have that $rv_1\neq 0.$
Therefore, by replacing $Q_1$ by $Q_1r$, $\theta_1(\cdot)$ by $\theta_1(\cdot)r$, $q$ by $r$ and $v_1$ by the partial isometry from the polar decomposition of $rv_1$, the relations \eqref{s1}-\eqref{s12} are still satisfied. This proves the claim.
\hfill$\square$

Let $i\in\{1,2\}$ and denote by $i+1$ the element in the set $\{1,2\}\setminus \{i\}.$
Since $P_i\prec_M Q_i$, by passing to relative commutants, we get that $A^{\Gamma_{i}}q_i\prec_ M P_{i+1}$. Then there exist projections $k_{i}\in A^{\Gamma_{i}}q_i$, $r_{i+1}\in P_{i+1}$, a $*$-homomorphism $\psi_i: A^{\Gamma_{i}}k_{i}\to r_{i+1}P_{i+1}r_{i+1}$ and a non-zero partial isometry 
$w_i\in r_{i+1}Mk_i$ such that 
\begin{equation}\label{c1}
\psi_i(x)w_i=w_ix, \text{    for all   } x\in A^{\Gamma_{i}}k_{i}.
\end{equation}
By restricting the projection $w_iw_i^*$ if necessary we can assume that 
\begin{equation}\label{c2}
E_{P_{i+1}}(w_iw_i^*)r_{i+1}\ge \beta_{i+1}r_{i+1}, \text{   for some positive number   } \beta_{i+1}.
\end{equation}
Denote $R_{i+1}=\psi_i(A^{\Gamma_{i}}k_{i}).$ 

{\bf Claim 2.}
$
R_1\prec_M^s A^{\Gamma_{2}}k_2 \text{   and   }
R_2\prec_M^s A^{\Gamma_{1}}k_1.
$
 
{\it Proof.} 
We prove only the second intertwining, since the first one follows in a similar way.
To the end, take a non-zero projection $b\in\mathcal N_{r_2Mr_2}(R_2)'\cap r_2Mr_2\subset R_2'\cap r_2P_2r_2 $ and define the $*$-homomorphism $\psi_b: A^{\Gamma_{1}}k_1\to R_2b$ by $\psi_b(x)=\psi_1(x)b, $ for all $x\in A^{\Gamma_{1}}k_1$. Note that $\psi_b(x)bw_1=bw_1x,$ for all $x\in A^{\Gamma_{1}}k_1$. We show that $bw_1$ is non-zero. If this is not the case, then $0=E_{P_2}(bw_1w_1^*)=bE_{P_2}(w_1w_1^*)$. Since $b\leq r_2$, relation \eqref{c2} implies that $b=0,$ false. Hence, $bw_1\neq 0.$ Thus, by replacing $bw_1$ by its partial isometry from its polar decomposition and by noticing that $\psi_b$ is a $*$-isomorphism, we get that $R_2b\prec_M A^{\Gamma_{1}}k_1.$ We use \cite[Lemma 2.4(2)]{DHI16} to conclude that $R_2\prec_M^s A^{\Gamma_{1}}k_1.$ 
\hfill$\square$

We continue with the following:

{\bf Claim 3.} $Q_1\vee A^{\Gamma_{1}}q_1\prec_M P_1\bar\otimes R_2$  \text{   and    }  $Q_2\vee A^{\Gamma_{2}}q_2\prec_M R_1\bar\otimes P_2.$

{\it Proof.} Due to symmetry we only need to show the first intertwining.
First we construct a $*$-isomorphism $\theta: Q_1k_1\to p_1P_1p_1$ and a non-zero partial isometry $v\in p_1Mk_1$ such that 
$$
\theta(x)v=vx, \text{   for all    } x\in Q_1k_1.
$$
Note that $\theta_1^{-1}:Q_1\to p_1P_1p_1$ is a $*$-isomorphism satisfying $\theta_1^{-1}(x)v_1^*=v_1^*x,$ for all $x\in Q_1.$ Since $Q_1$ is a factor and $k_1\in Q_1'\cap q_1(A\rtimes\Gamma_{1})q_1$, we can define the $*$-homomorphism $\theta: Q_1k_1\to p_1P_1p_1$, by letting $\theta(xk_1)=\theta_1^{-1}(x),$ for all $x\in Q_1.$ Note that $\theta(xk_1)v^*_1k_1=v^*_1k_1xk_1,$ for all $x\in Q_1.$
Let $v$ be the partial isometry obtained from the polar decomposition of $v_1^*k_1$. Note that $v\neq 0$ by using \eqref{s12}. Therefore, $\theta(x)v=vx$, for all $x\in Q_1k_1.$

Notice that $vv^*\in (P_1'\cap M)p_1$ and $P_1'\cap M= P_2$, so there exists a non-zero projection $\tilde p_2\in P_2$ such that $vv^*=p_1\otimes \tilde p_2$.
Without loss of generality, we can assume that 
\begin{equation}\label{s4}
\tilde p_2 \leq r_2 \text{    or    } r_2 \leq \tilde p_2. 
\end{equation}
Indeed, since $P_2$ is a II$_1$ factor, there exists a unitary $u\in\mathcal U(P_2)$ such that $\tilde p_2\leq u r_2 u^*$ or $ur_2u^*\leq \tilde p_2.$ 
By replacing $\psi_1(\cdot)$ by $u\psi_1(\cdot) u^*$, $r_2$ by $ur_2u^*$, $R_2$ by $uR_2u^*$ and $w_1$ by $uw_1$, relations \eqref{c1} and \eqref{c2} still hold. Therefore we can assume \eqref{s4} to be true.

We suppose by contradiction that $Q_1k_1\vee A^{\Gamma_{1}}k_1\nprec_M P_1\bar\otimes R_2.$ Then there exist two sequences of unitaries $(u_n)_n\subset\mathcal U(Q_1)$ and $(v_n)_n\subset \mathcal U(A^{\Gamma_{1}})$ such that 
$$
\|E_{P_1\bar\otimes R_2}(x(u_nk_1)(v_nk_1)y)\|_2\to 0, \text{ for all  }  x,y\in M.
$$
By taking $x=v$, we get that 
\begin{equation}\label{e}
\|E_{P_1\bar\otimes R_2}(\theta(u_nk_1)v v_ny)\|_2=
\|E_{P_1\bar\otimes R_2}(v v_ny)\|_2\to 0, \text{ for all  }  y\in M
\end{equation}
We now argue that 
\begin{equation}\label{extra}
E_{P_1\bar\otimes R_2}(vv^* w_1w_1^*)=0.
\end{equation}

Take $y=au_g^*w_1^*$ in \eqref{e}, for some $a\in A$ and $g\in\Gamma$. Since $A^{\Gamma_{1}}$ is normalized by $u_g$, we obtain
$$
\begin{array}{rcl}
\|E_{P_1\bar\otimes R_2}(v v_nau_g^*w_1^*)\|_2&=&
\|E_{P_1\bar\otimes R_2}(v au_g^*(\sigma_g(v_n)w_1^*))\|_2\\
&=&\|E_{P_1\bar\otimes R_2}(v au_g^*w_1^*\psi_1(\sigma_g(v_n)k_1))\|_2\\
&=&\|E_{P_1\bar\otimes R_2}(v au_g^*w_1^*)\|_2\\
\end{array}
$$
tends to $0$ as $n\to\infty$. This implies that $\|E_{P_1\bar\otimes R_2}(v au_g^*w_1^*)\|_2=0$, for all $a\in A$ and $g\in\Gamma$. This proves \eqref{extra}, which gives us that
$p_1E_{P_1\bar\otimes R_2}(\tilde p_2w_1w_1^*)=0$. We obtain $p_1\otimes E_{R_2}(\tilde p_2w_1w_1^*)=0$ and therefore,
$p_1\otimes E_ {R_2}(\tilde p_2E_{P_2}(w_1w_1^*))=0$, since $\tilde p_2\in P_2$. Hence, $\tau(r_2\tilde p_2E_{P_2}(w_1w_1^*))=0$. Using that $E_{P_2}(w_1w_1^*)r_2\ge \beta_2r_2$, we get that $\tilde p_2E_{P_2}(w_1w_1^*)r_2\tilde p_2\ge \beta_2r_2\tilde p_2$. Altogether it implies that $\beta_2\tau(\tilde p_2r_2)=0$, which contradicts \eqref{s4} since $\beta_2$ is non-zero. Thus, $Q_1k_1\vee A^{\Gamma_{1}}k_1\prec_M P_1\bar\otimes R_2$, ending this way the proof of the claim.
\hfill$\square$

Finally, we obtain the following:

{\bf Claim 4.} $A\prec_M^s A^{\Gamma_{1}}\vee A^{\Gamma_{2}}.$

{\it Proof.} 
Denote $D_i=R_i\oplus \mathbb C (1-r_i),$ for $i\in\{1,2\}.$ Note that Claim 1 together with relation \eqref{s1} show that $Q_i\vee A^{\Gamma_i}q_1\subset q_i(A\rtimes\Gamma_i)q_i$ has finite index, for $i\in\{1,2\}.$
Combining Claim 3 and Lemma \ref{L: center}(2), we obtain that $A\prec_M P_1\bar\otimes D_2$ and $A\prec_M D_1\bar\otimes P_2.$ Note that \cite[Lemma 2.4(2)]{DHI16} together with Proposition \ref{L: PV} imply that $A\prec^s_M D_1\bar\otimes D_2.$ Notice also that Claim 2 shows that $D_i\prec_M^s A^{\Gamma_{1}}\vee A^{\Gamma_{2}}$, for any $i\in\{1,2\}.$ Since $A^{\Gamma_{1}}\vee A^{\Gamma_{2}}$ is regular and the algebra $D_2$ is abelian, we can apply Lemma \ref{L: joint} and obtain that $D_1\bar\otimes D_2\prec_M^s A^{\Gamma_{1}}\vee A^{\Gamma_{2}}.$ By applying \cite[Lemma 3.7]{Va08}, we deduce that $A\prec_M^s A^{\Gamma_{1}}\vee A^{\Gamma_{2}}$. 
\hfill$\square$

The moreover part follows from the first part of Theorem \ref{Th:split}. 
In particular, we can represent $A^{\Gamma_2}=L^\infty(X_1,\mu_1)$ and $A^{\Gamma_2}=L^\infty(X_2,\mu)$ for some standard probability spaces $(X_1,\mu_1)$ and $(X_2,\mu_2)$, respectively. Hence, for any $i\in\{1,2\}$ there exists a pmp action $\Gamma_{i}\car (X_i,\mu_i)$ such that $\Gamma\car X$ is isomorphic to $\Gamma_{1}\times\Gamma_{2}\car X_1\times X_2. $
\hfill$\blacksquare$

\section{Proofs of the main results}

We start this section by presenting another tool needed for the proof of Theorem \ref{A} and we will conclude by proving the main results mentioned in the introduction. 

\begin{proposition}\label{start}
Let $\Gamma=\Gamma_1\times\dots \times\Gamma_n$ be a product of $n\ge 2$ groups that belong to the class $\mathcal C$. Let $\Gamma\overset{}{\car} (X,\mu)$ be a free ergodic pmp action and 
denote $M=L^\infty (X)\rtimes\Gamma.$\\
Suppose that $M= P_1\bar\otimes P_2 $, for some II$_1$ factors $P_1$ and $P_2.$

Then there exists a partition $T_1\sqcup T_2=\{1,\dots ,n \}$ into non-empty sets such that $P_i\prec^s_M L^\infty(X)\rtimes\Gamma_{T_i}$, for $i\in\{1,2\}.$
\end{proposition}

{\it Proof.} Denote $A=L^\infty(X).$ Let $T_i$ be the minimal subset of $\{1,\dots,n\}$ with the property that $P_i \prec^s_M A\rtimes\Gamma_{T_i}$ for all $i\in \{1,2\}.$ Notice that $T_i$ is non-empty since any corner of a II$_1$ factor is non-abelian. We want to show that $T_1 \sqcup T_2=\{1,\dots,n\}.$ Note that since $P_i$ is regular and $M$ is a factor, \cite[Lemma 2.4(2)]{DHI16} implies that $P_i\prec_M A\rtimes{\Gamma_S}$ if and only if $P_i\prec_M^s A\rtimes {\Gamma_S}$ for any subset $S\subset \{1,\dots,n\}$.

First we notice that $\{1,\dots,n\}=T_1\cup T_2$. Indeed, by applying \cite[Lemma 2.3]{BV12} we get that $M\prec A\rtimes \Gamma_{T_1\cup T_2}$. This shows that $\{1,\dots,n\}=T_1\cup T_2$, since $\Gamma_t$ is an infinite group, for any $t\in\{1,\dots,n\}.$ We will finish the proof by proving the following claim.

{\bf Claim.} $T_1\cap T_2$ is empty. 

{\it Proof.} For any $i\in\{1,2\},$
\cite[Proposition 2.4]{CKP14} implies that there exist non-zero projections $p_i\in P_i$, $q_i\in A\rtimes\Gamma_{T_i}$, a subalgebra $Q_i\subset q_i(A\rtimes\Gamma_{T_i})q_i$, a partial isometry $v_i\in q_iMp_i$ and an onto $*$-
 $\theta_i: p_iP_ip_i\to Q_i$ such that $\theta_i(x)v_i=v_ix, \text{   for all   } x\in p_iP_ip_i, $ and
$
Q_i\vee (Q_i'\cap q_i(A\rtimes\Gamma_{T_i})q_i)\subset q_i(A\rtimes\Gamma_{T_i})q_i \text{  has finite index. }
$
Moreover, the support projection of $E_{A\rtimes\Gamma_{T_i}}(v_iv_i^*)$ can be assumed to equal $q_i.$ Denote $S_i:=Q_i'\cap q_i(A\rtimes\Gamma_{T_i})q_i$.

Assume by contradiction that there exist a non-zero projection $z\in\mathcal N_{q_1(A\rtimes\Gamma_{T_1})q_1}(S_1)'\cap q_1(A\rtimes\Gamma_{T_1})q_1$ and an index $j\in T_1$ such that $S_1z_0$ is non-amenable relative to $A\rtimes \Gamma_{T_{1}\setminus\{j\}},$ for all $z_0\in (Q_1\vee S_1)'\cap z(A\rtimes\Gamma_{T_1})z$. Note that $z,z_0\in S_1$ and that the inclusion $z_0(Q_1\vee S_1)z_0\subset z_0(A\rtimes\Gamma_{T_1})z_0$
has finite index by Lemma \ref{PP}. Therefore,
Corollary \ref{all} implies that $Q_1z_0\prec_{A\rtimes{{\Gamma_{T_1}}}} A\rtimes \Gamma_{T_{1}\setminus\{j\}}.$ By applying \cite[Lemma 2.4(3)]{DHI16}, we get that $Q_1z\prec^s_{A\rtimes{{\Gamma_{T_1}}}} A\rtimes \Gamma_{T_{1}\setminus\{j\}}.$ It is easy to see that the moreover part of the previous paragraph shows that $P_1\prec_M Q_1z.$ Hence, by applying \cite[Lemma 3.7]{Va08} we get that $P_1\prec_M A\rtimes\Gamma_{T_{1}\setminus\{j\}},$ which contradicts the minimality of $T_1.$

Therefore, for any $j\in T_1$ and $z\in\mathcal N_{q_1(A\rtimes\Gamma_{T_1})q_1}(S_1)'\cap q_1(A\rtimes\Gamma_{T_1})q_1$, there exists a non-zero projection $z_0\in (Q_1\vee S_1)'\cap z(A\rtimes\Gamma_{T_1})z$ such that $S_1z_0$ is amenable relative to $A\rtimes \Gamma_{T_{1}\setminus\{j\}}.$ \cite[Lemma 2.6(2)]{DHI16} shows that we can assume $z_0\in \mathcal N_{z(A\rtimes\Gamma_{T_1})z}(S_1z)'\cap z(A\rtimes\Gamma_{T_1})z$. By applying \cite[Proposition 2.7]{PV11} finitely many times, we obtain that there exists a non-zero projection $z_1\in\mathcal N_{q_1(A\rtimes\Gamma_{T_1})q_1}(S_1)'\cap q_1(A\rtimes\Gamma_{T_1})q_1$ such that $S_1z_1$ is amenable. In a similar way, there exists a non-zero projection $z_2\in\mathcal N_{q_2(A\rtimes\Gamma_{T_2})q_2}(S_2)'\cap q_2(A\rtimes\Gamma_{T_2})q_2\subset Q_2$ such that $S_2z_2$ is amenable.

Since $P_1\prec_M Q_1z_1$ and $P_2\prec_M Q_2z_2$, \cite[Lemma 3.5]{Va08} implies that $S_1z_1\prec_M P_2$ and $S_2z_2\prec_M P_1.$ By proceeding as in the proof of Theorem \ref{AA} (Claims 3 and 4), we obtain that there exist amenable subalgebras $D_1\subset P_1$ and $D_2\subset P_2$ such that $L(\Gamma_{T_1}\cap \Gamma_{T_2})\prec_M D_1\bar\otimes D_2.$ This implies that  $\Gamma_{T_1}\cap \Gamma_{T_2}$ is amenable, hence $T_1\cap T_2$ is empty.
\hfill$\blacksquare$
 
 \begin{remark}
 We provide the following shorter argument for proving Proposition \ref{start} if the groups $\Gamma_i$'s are only weakly amenable, bi-exact groups or free products.  We only need to prove the claim. First, for any $t\in \{1,\dots,n\}$, denote $\hat t=\{1,\dots,n\}\setminus\{t\}$. Suppose by contradiction that there exists an element $t\in T_1\cap T_2.$ Then $P_1\nprec_M A\rtimes \Gamma_{T_1\setminus\{t\}}$ based on the minimality of $T_1$. This is equivalent to $P_1\nprec_M A\rtimes \Gamma_{\hat t}$ using \cite[Lemma 2.8(2)]{DHI16}. If $\Gamma_t$ is a free product, by applying Lemma \ref{L: amalgam}, we must have $P_2\prec_M A\rtimes\Gamma_{\hat t}$. Using \cite[Lemma 2.8(2)]{DHI16} this shows that $P_2\prec_M A\rtimes\Gamma_{T_2\setminus\{t\}}$, which contradicts the minimality of $T_2$. On the other hand, if $\Gamma_t\in \mathcal C_{rss}$, by applying Lemma \ref{L: rss} we obtain that $P_2$ is amenable relative to $A\rtimes\Gamma_{\hat t}$, which implies that $P_2\prec A \rtimes\Gamma_{T_2\setminus \{t\}}$ or $M$ is amenable relative to $A\rtimes\Gamma_{\hat t}.$ The former contradicts the minimality of $T_2$, while the last one contradicts the non-amenability of $\Gamma_t$ by \cite[Proposition 2.4]{OP07}.  
 \end{remark}

The proof of Theorem \ref{A} follows by combining Proposition \ref{start} with Theorem \ref{AA}. 
Corollary \ref{A2} is obtained directly from Theorem \ref{A}. We continue now with the proofs of Corollary \ref{C} and Theorem \ref{UPFgeneral}.

{\it Proof of Corollary \ref{C}.} 
By applying Theorem \ref{A} finitely many times, we can find an integer $1\leq k \leq n$, a partition $S_1\sqcup ... \sqcup S_k =\{1,...,n\}$ and pmp actions $\Gamma_{S_{i}}\car (X_i,\mu_i)$ such that 
$\Gamma\car X \text{  is isomorphic to } \Gamma_{S_1}\times ...\times \Gamma_{S_k}\car X_1\times ...\times X_k,
$
and $M_i=L^\infty(X_i)\rtimes\Gamma_{S_i}$ is prime for all $i\in\{1,...,k\}$. Note that the following holds:

{\bf Claim.} If $\Gamma\car X$ is isomorphic to $\Gamma_{T_1}\times\Gamma_{T_2}\car Y_1\times Y_2$ for a partition $T_1\sqcup T_2=\{1,...,n\}$, then there exists a partition $J_1\sqcup J_2=\{1,...,k\}$ such that
$T_1=\sqcup_{i\in J_1} S_i$ and $T_2=\sqcup_{i\in J_2} S_i.$

{\it Proof.} First note that it is enough to show that if $S_i\cap T_j\neq \emptyset $, for some $i\in\{1,...,k\}$ and $j\in\{1,2\}$, then $S_i\subset T_j.$ To this end, take $i$ and $j$ as before. The assumption implies that the actions $\Gamma_{S_i}\car X$ and $\Gamma_{S_i\cap T_1}\times\Gamma_{S_i\cap T_2}\car Y_1\times Y_2$ are isomorphic. This shows that 
$$
M_i=L^\infty(X_i)\rtimes\Gamma_{S_i}=L^\infty(X)^{\Gamma_{S_i^c\cap T_2}}\rtimes\Gamma_{S_i\cap T_1}\bar\otimes L^\infty(X)^{\Gamma_{S_i^c\cap T_1}}\rtimes\Gamma_{S_i\cap T_2},
$$
where we denote $S_i^c:=\{1,...,n\}\setminus S_i$, for any $1\leq i\leq k.$
Since the algebra $M_i$ is prime, we must have $S_i\cap T_j=S_i,$ which implies that $S_i\subset T_j.$
\hfill$\square$

The Claim shows that the partition $S_1\sqcup ... \sqcup S_k =\{1,...,n\}$ is unique up to a permutation of the sets. We continue now by proving the moreover part. 

(1) Let $M=P_1\bar\otimes P_2$ for some II$_1$ factors $P_1$ and $P_2$. If we apply Theorem \ref{A} we obtain a partition $T_1\sqcup T_2=\{1,...,n\}$, a decomposition $M=P_1^t\bar\otimes P_2^{1/t}$ for some positive $t$, and a unitary $u\in M$ such that 
$$
P_1^{t}=u(L^\infty(X)^{\Gamma_{T_2}}\rtimes\Gamma_{T_1})u^* \text{  and  }
P_2^{1/t}=u(L^\infty(X)^{\Gamma_{T_1}}\rtimes\Gamma_{T_2})u^*.
$$ 
The Claim shows that there exists a partition $J_1\sqcup J_2=\{1,...,k\}$ such that
$T_1=\sqcup_{i\in J_1} S_i$ and $T_2=\sqcup_{i\in J_2} S_i.$ This exactly implies the conclusion.

(2) Assume that $M=P_1\bar\otimes \dots \bar\otimes P_m$, for some $m\ge k$. Part (1) combined with induction implies that $m\leq k$ and that there exist a partition $J_1\sqcup\dots \sqcup J_m=\{1,\dots,k\}$, a decomposition $M=P_1^{t_1}\bar\otimes\dots\bar\otimes P_m^{t_m}$ with $t_1\dots t_m=1$, and a unitary $u\in M$ such that 
$P_s^{t_j}=u(\bar\otimes_{j\in J_s} M_j)u^*$, for any $s\in\{1,...,m\}.$ Therefore $m=k$ and the conclusion holds.

(3) For proving this last part, we proceed as in (2). Since each $P_j$ is prime, $J_s$ has only one element. This shows once again that $m=k$ and the conclusion holds.
\hfill$\blacksquare$

{\it Proof of Theorem \ref{UPFgeneral}.} 
(1) Denote $X=X_1\times ...\times X_k$. Let $M=P_1\bar\otimes P_2$ for some II$_1$ factors $P_1$ and $P_2$. By applying Proposition \ref{start} and Theorem \ref{Th:split}(2), we obtain that there exists a partition $I_1\sqcup I_2=\{1,...,n\}$ such that $P_1$ is stably isomorphic to $L^\infty(X)^{\Gamma_{I_2}}\rtimes\Gamma_{I_1}=\bar\otimes_{i\in I_1}M_i$ and $P_2$ is stably isomorphic to $L^\infty(X)^{\Gamma_{I_1}}\rtimes\Gamma_{I_2}=\bar\otimes_{i\in I_2}M_i.$

(2) \& (3) Assume that $M=P_1\bar\otimes\dots\bar\otimes P_m$ for some integer $m$ and II$_1$ factors $P_1,\dots ,P_m.$ Part (1) combined with induction implies that $m\leq k$ and that there exists a partition $J_1\sqcup\dots \sqcup J_m=\{1,\dots,k\}$ such that $P_s$ is stably isomorphic to $\bar\otimes_{j\in J_s}M_j$, for any $1\leq s\leq m$. If $m\ge k$ or if each $P_s$ is prime, we obtain that $m=k$ and each $J_s$ has only one element. 

For proving the moreover part, assume that $\Gamma_i\car (X_i,\mu_i)$ is strongly ergodic for any $1\leq i\leq k.$ Note that \cite[Theorem C]{CSU13} combined with \cite[Examples 1.4,1.5]{CSU13} imply that the class $\mathcal C$ is contained in the class of non-inner amenable groups. In combination with \cite{Ch82}, we obtain that $M$ does not have property Gamma. 

Let $M=P_1\bar\otimes P_2$ for some II$_1$ factors $P_1$ and $P_2$. We will show only part (1) from the moreover part since part (2) and (3) can be deduced as in Corollary \ref{C}. The previous paragraph implies that $P_1$ and $P_2$ does not have property Gamma. As before, we can apply Proposition \ref{start} and obtain a partition $I_1\sqcup I_2=\{1,...,n\}$ such that $P_j\prec_M L^\infty(X)\rtimes I_j$, for any $j\in\{1,2\}$. Applying \cite[Proposition 6.3]{Ho15}, we obtain that $P_j\prec_M \bar\otimes_{i\in I_j} M_i$, for any $j\in\{1,2\}.$ By proceeding as in the proof of Theorem \ref{AA}, we obtain that there exist a unitary $u\in M$ and a decomposition $M=P_1^t\bar\otimes P_2^{1/t}$, for some $t>0$, such that $P_1^t= u(\bar\otimes_{i\in I_1} M_i)u^*$ and 
$P_2^{1/t}= u(\bar\otimes_{i\in I_2} M_i)u^*$. This ends the proof.
\hfill$\blacksquare$


\begin{thebibliography}{ABC99}


\bibitem[Bo12]{Bo12} R. Boutonnet: {\it On solid ergodicity for Gaussian actions}, J. Funct. Anal. {\bf 263} (2012), no. 4, 1040-1063.

\bibitem[BV12]{BV12} M. Berbec, S. Vaes: {\it W$^*$-superrigidity for group von Neumann algebras of left-right wreath products}, Proceedings of the London Mathematical Society (3) {\bf 108} (2014), no. 5, 1116-1152.



\bibitem[CdSS15]{CdSS15} I. Chifan, R. de Santiago, T. Sinclair: {\it W$^*$-rigidity for the von Neumann algebras of products of hyperbolic groups}, Geom. Funct. Anal. {\bf 26} (2016), no. 1, 136-159.

\bibitem[CdSS17]{CdSS17} I. Chifan, R. de Santiago,  W. Sucpikarnon: {\it Tensor product decompositions of II$_1$ factors arising from extensions of amalgamated free product groups}, Communications in Mathematical Physics {\bf 364} (2018), 1163-1194.


\bibitem[Ch82]{Ch82} M. Choda: {\it Inner amenability and fullness}, Proc. Amer. Math. Soc. {\bf 86} (1982), 663-666.

 \bibitem[CH08]{CH08} I. Chifan, C. Houdayer: {\it Bass-Serre rigidity results in von Neumann algebras}, Duke Math. J. {\bf 153} (2010), no. 1, 23-54.

\bibitem[CI08]{CI08} I. Chifan, A. Ioana: {\it Ergodic subequivalence relations induced by a Bernoulli action}, Geom. Funct. Anal. {\bf 20} (2010), no. 1, 53-67.





 \bibitem[CIK13]{CIK13} I. Chifan, Y. Kida, A. Ioana: {\it W$^*$-superrigidity for arbitrary actions of central quotients of braid groups},
  Math. Ann. {\bf 361} (2015), no. 3-4, 563-582.

\bibitem[CKP14]{CKP14} I. Chifan, Y. Kida, S. Pant: {\it Primeness results for von neumann algebras associated with surface braid groups},  Int. Math. Res. Not. (2015), article id:rnv271, 42pp.

\bibitem[Co76]{Co76} A. Connes, Classification of injective factors, Ann. Math. {\bf 104} (1976) 73-115.

\bibitem[CPS11]{CPS11} I. Chifan, S. Popa, J. O. Sizemore, {\it Some OE- and W$^*$
-rigidity results for actions by wreath
product groups}, J. Funct. Anal. {\bf 263} (2012), no. 11, 3422-3448.

 \bibitem[CS11]{CS11} I. Chifan, T. Sinclair: {\it On the structural theory of II$_1$ factors of negatively curved groups}, Ann. Sci. \'{E}c. Norm. Sup\'{e}r. (4) 46 (2013), no. 1, 1-33 (2013).
 
 \bibitem[CSU13]{CSU13}  I. Chifan, T. Sinclair and B. Udrea: {\it Inner amenability for groups and central sequences in factors}, Ergodic Theory Dynam. Systems {\bf 36} (2016), 1106-1029.




\bibitem[De19]{De19} T. Deprez: {Ozawa’s class $\mathcal S$ for locally compact groups and unique prime
factorization of group von Neumann algebras}, arXiv:1904.02090.


\bibitem[DHI16]{DHI16} D. Drimbe, D. Hoff, A. Ioana: {\it Prime II$_1$ factors arising from irreducible lattices in product of rank one simple Lie groups}, to appear in J. Reine. Angew. Math, arXiv: 1611.02209 (2016).

\bibitem[DI12]{DI12} Y. Dabrowski, A. Ioana: {\it Unbounded derivations, free dilations, and indecomposability results for II$_1$ factors}, Trans. Amer. Math. Soc. {\bf 368} (2016), no. 7, 4525-4560.







\bibitem[Dr19]{Dr19} D. Drimbe: {\it Orbit equivalence rigidity for product actions}, to appear in Comm. Math. Physics., arXiv: 1905.07642.

\bibitem[dSP18]{dSP18} R. de Santiago, S. Pant: {\it Classification of tensor decompositions of II$_1$ factors associated with poly-hyperbolic groups}, arXiv: 1802.09083(2018).

\bibitem[Ge95]{Ge95} L. Ge: {\it On maximal injective subalgebras of factors}, Adv. Math. {\bf 118} (1996), no. 1, 34-70.

\bibitem[Ge96]{Ge96} L. Ge: {\it Applications of free entropy to finite von Neumann algebras}, II, Ann. of Math. (2) {\bf 147} (1998), no. 1, 143-157.


\bibitem[HI15]{HI15} C. Houdayer, Y. Isono: {\it Unique prime factorization and bicentralizer problem for a class
of type III factors},  Adv. Math. {\bf 305} (2017), 402-455.

\bibitem[Ho15]{Ho15} D. Hoff: {\it Von Neumann algebras of equivalence relations with nontrivial one-cohomology}, J. Funct. Anal. {\bf 270} (2016), no. 4, 1501-1536.

\bibitem[HPV10]{HPV10} C. Houdayer, S. Popa, S. Vaes: {\it A class of groups for which every action is W$^*$-superrigid}, Groups Geom. Dyn. {\bf 7} (2013), no. 3, 577-590.

\bibitem[HV12]{HV12} C. Houdayer, S. Vaes: {\it Type III factors with unique Cartan decomposition},  J. Math\'ematiques Pures et Appliqu\'ees {\bf 100} (2013), 564-590.





\bibitem [Io06]{Io06} A. Ioana: {\it Rigidity results for wreath product of II$_1$ factors}, J. Funct. Anal. {\bf 252} (2007), 763-791.

\bibitem [Io10]{Io10} A. Ioana: {\it W$^*$-superrigidity for Bernoulli actions of property (T) groups}, J. Amer. Math. Soc. {\bf 24} (2011),
1175-1226.

\bibitem[Io11]{Io11} A. Ioana: {\it Uniqueness of the group measure space decomposition for Popa's $\mathcal H\mathcal T$ factors}, Geom. Funct. Anal. {\bf 22} (2012), no. 3, 699-732.

\bibitem [Io12a]{Io12a} A. Ioana: {\it Classification and rigidity for von Neumann algebras}, European Congress of Mathematics, EMS (2013), 601-625.

\bibitem[Io12]{Io12} A. Ioana: {\it Cartan subalgebras of amalgamated free product of II$_1$ factors}, With an appendix by Adrian Ioana and Stefaan Vaes, Ann. Sci. \'Ec. Norm. Sup\'er. (4) {\bf 48} (2015), no. 1, 71-130.

\bibitem [Io14]{Io14} A. Ioana: {\it Strong ergodicity, property (T), and orbit equivalence rigidity for translation actions}, J. Reine. Angew. Math. {\bf 733} (2017), 203-250.

\bibitem[Io17]{Io17} A. Ioana: {\it Rigidity for von Neumann algebras}, submitted to Proceedings ICM 2018, arXiv: 1712.00151.

\bibitem[IPV10]{IPV10} A. Ioana, S. Popa, S. Vaes: {\it A class of superrigid group von Neumann algebras}, Ann. of Math. (2) {\bf 178} (2013), no. 1, 231-286.

 \bibitem[Is14]{Is14} Y. Isono: {\it Some prime factorization results for free quantum group factors}, J. Reine. Angew. Math. {\bf 722}, (2017), 215-250. 
 
 
 \bibitem[Is16]{Is16} Y. Isono: {\it On fundamental groups of tensor product II$_1$ factors}, to appear in J. Inst. Math. Jussieu, arXiv:1608.06426.
 


\bibitem[Jo81]{Jo81} V.F.R. Jones, {\it Index for subfactors}, Invent. Math. {\bf 72} (1983), 1-25.


\bibitem[KV15]{KV15} A. Krogager, S. Vaes: {\it A class of II$_1$ factors with exactly two crossed product decompositions}, J. Math\'ematiques Pures et Appliqu\'ees, {\bf 108} (2017), 88-110.


\bibitem[MvN36]{MvN36} F.J. Murray, J. von Neumann, On rings of operators, Ann. Math. {\bf 37} (1936), 116-229.

\bibitem[MvN43]{MvN43} F.J. Murray, J. von Neumann, Rings of operators IV, Ann. Math. {\bf 44} (1943), 716-808.






 \bibitem[OP03]{OP03} N. Ozawa, S. Popa: {\it Some prime factorization results for type II$_1$ factors},  Invent. Math. {\bf 156} (2004), no. 2, 223-234.

\bibitem[OP07]{OP07} N. Ozawa, S. Popa: {\it On a class of II$_1$ factors with at most one Cartan subalgebra}, Ann. of Math. (2) {\bf 172} (2010), no. 1, 713-749.


\bibitem [OW80]{OW80} D. Ornstein, B. Weiss: {\it Ergodic theory of amenable group actions. I. The Rohlin lemma}, Bull. Amer. Math. Soc. (N.S.), 2(1):161-164, 1980. 


\bibitem[Oz03]{Oz03} N. Ozawa: {\it Solid von Neumann algebras}, Acta Math. {\bf 192} (2004), no. 1, 111-117.



\bibitem[Oz04]{Oz04} N. Ozawa: {\it A Kurosh type theorem for type II$_1$ factors}, Int. Math. Res. Not. (2006), Art. ID 97560, 21
pp.


 \bibitem[Pe06]{Pe06} J. Peterson: {\it L$^2$-rigidity in von Neumann algebras}, Invent. Math. {\bf 175} (2009), no. 2, 417-433.

\bibitem[Po81]{Po81} S. Popa: {\it On a problem of R. V. Kadison on maximal abelian $*$-subalgebras in factors}, Invent. Math {\bf 65} (1981/82), no. 2, 269-281. MR 641132 (83g:46056).

\bibitem[Po83]{Po83} S. Popa: {\it Orthogonal pairs of $*$-subalgebras in finite von Neumann algebras}, J. Operator Theory {\bf 9} (1983), no. 2, 253-268.

\bibitem [Po01]{Po01} S. Popa: {\it On a class of type II$_1$ factors with Betti numbers invariants}, Ann. of Math. {\bf 163} (2006), 809-899.

\bibitem[Po03]{Po03} S. Popa: {\it Strong rigidity of II$_1$ factors arising from malleable actions of {\it w}-rigid groups. I.}, Invent. Math. {\bf 165} (2006), no. 2, 369-408.


 \bibitem[Po06a]{Po06a} S. Popa: {\it On the superrigidity of malleable actions with spectral gap}, J. Amer. Math. Soc. {\bf 21} (2008),
981-1000.

\bibitem[Po06b]{Po06b} S. Popa: {\it On Ozawa's property for free group factors}, Int. Math. Res. Not. 
2007, no. 11, Art. ID rnm036, 10 pp.

\bibitem [Po07]{Po07} S. Popa: {\it Deformation and rigidity for group actions and von Neumann algebras}, In Proceedings of the ICM (Madrid, 2006), Vol. I, European Mathematical Society Publishining House, 2007, 445-477.

\bibitem[PP86]{PP86} M. Pimsner, S. Popa: {\it Entropy and index for subfactors}, Ann. Sci. \'Ecole Norm. Sup. {\bf 19} (1986), 57=106.




\bibitem[PS03]{PS03} S. Popa, D. Shlyakhtenko: {\it Cartan subalgebras and bimodules decompositions of II$_1$ factors}, Math. Scand. {\bf 92} (2003), 93-102.

\bibitem[PV11]{PV11} S. Popa, S. Vaes: {\it Unique Cartan decomposition for II$_1$ factors arising from arbitrary actions of free groups}, Acta Math. {\bf 212} (2014), no. 1, 141-198.
 
\bibitem[PV12]{PV12} S. Popa, S. Vaes: {\it Unique Cartan decomposition for II$_1$ factors arising from arbitrary actions of hyperbolic groups}, J. Reine Angew. Math. {\bf 694} (2014), 215-239.

\bibitem[Si10]{Si10} T. Sinclair: {\it Strong solidity of group factors from lattices in SO(n, 1) and SU(n, 1)}, J. Funct. Anal. {\bf 260} (2011), no. 11, 3209-3221. MR 2776567 (2012e:22009).



 \bibitem[SW11]{SW11}  O. Sizemore, A. Winchester: {\it Unique prime decomposition results for factors coming
from wreath product groups}, Pacific J. Math. {\bf 265} (2013),  221-232. 

\bibitem[Ta01]{Ta01} M. Takesaki: {\it Theory of Operator Algebras I}, ser. Encyclopaedia of Mathematical Sciences. Springer, {\bf 124},
2001. xix+415 pp.




\bibitem[Va08]{Va08} S. Vaes: {\it Explicit computations of all finite index bimodules for a family of II$_1$ factors}, Ann. Sci. \'{E}c. Norm. Sup\'{e}r. (4) {\bf 41} (2008), no. 5, 743-788.



\bibitem [Va10a]{Va10a} S. Vaes: {\it Rigidity for von Neumann algebras and their invariants}, Proceedings of the ICM (Hyderabad, India, 2010), Vol. III, Hindustan Book Agency (2010), 1624-1650.

\bibitem[Va10b]{Va10b} S. Vaes: {\it One-cohomology and the uniqueness of the group measure space decomposition of a II$_1$ factor}, Math. Ann. {\bf 355} (2013), no. 2, 661-696.

\bibitem[Va13]{Va13} S. Vaes: {\it Normalizers inside amalgamated free product von Neumann algebras}. Publ. Res. Inst. Math.
Sci. {\bf 50} (2014), 695-721.


\end{thebibliography}
\end{document}